\documentclass{amsart}
\usepackage{amsfonts}

\newtheorem{theorem}{Theorem}
\theoremstyle{plain}

\newtheorem{corollary}{Corollary}

\newtheorem{proposition}{Proposition}
\newtheorem{remark}{Remark}

\numberwithin{equation}{section}
\input{tcilatex}

\begin{document}
\title[Cauchy-Bunyakovsky-Schwarz  type Inequalities]{Some
Cauchy-Bunyakovsky-Schwarz Type Inequalities for Sequences of Operators in
Hilbert Spaces}
\author{Sever S. Dragomir}
\address{School of Computer Science and Mathematics\\
Victoria University of Technology\\
PO Box 14428, MCMC 8001\\
Victoria, Australia.}
\email{sever@matilda.vu.edu.au}
\urladdr{http://rgmia.vu.edu.au/SSDragomirWeb.html}

\begin{abstract}
Some inequalities of Cauchy-Bunyakovsky-Schwarz  type for sequences of
bounded linear operators in Hilbert spaces and some applications are given.
\end{abstract}

\maketitle

\section{Introduction}

Let $\left( H;\left( \cdot ,\cdot \right) \right) $ be a real or complex
Hilbert space and \thinspace $B\left( H\right) $ the Banach algebra of all
bounded linear operators that map $H$ into $H.$

We recall that a self-adjoint operator $A\in B\left( H\right) $ is positive
in $B\left( H\right) $ iff $\left( Ax,x\right) \geq 0$ for any $x\in H.$ The
binary relation $A\geq B$ iff $A-B$ is a positive self-adjoint operator, is
an \textit{order relation }\ on $B\left( H\right) .$ We remark that for any $%
A\in B\left( H\right) $ the operators $U:=AA^{\ast }$ and $V:=A^{\ast }A$
are positive self adjoint operators on $H$ and $\left\Vert U\right\Vert
=\left\Vert V\right\Vert =\left\Vert A\right\Vert ^{2}.$

In \cite{SSD}, the author has proved the following inequality of
Cauchy-Bunyakovsky-Schwarz type in the order of $B\left( H\right) .$

\begin{theorem}
\label{t1.1}Let $A_{1},\dots ,A_{n}\in B\left( H\right) $ and $z_{1},\dots
,z_{n}\in \mathbb{K}$ $\left( \mathbb{R},\mathbb{C}\right) .$ Then the
following inequality holds:%
\begin{equation}
\sum_{i=1}^{n}\left\vert z_{i}\right\vert ^{2}\sum_{i=1}^{n}A_{i}A_{i}^{\ast
}\geq \left( \sum_{i=1}^{n}z_{i}A_{i}\right) \left( \sum_{i=1}^{n}\overline{%
z_{i}}A_{i}^{\ast }\right) \geq 0.  \label{1.1}
\end{equation}
\end{theorem}

\begin{proof}
For the sake of completeness, we give here a simple proof of this inequality.

For any $i,j\in \left\{ 1,\dots ,n\right\} $ one has in the order of $%
B\left( H\right) :$%
\begin{equation*}
\left( \overline{z_{i}}A_{j}-\overline{z_{j}}A_{i}\right) \left( \overline{%
z_{i}}A_{j}-\overline{z_{j}}A_{i}\right) ^{\ast }\geq 0,
\end{equation*}%
that is,%
\begin{equation*}
\left( \overline{z_{i}}A_{j}-\overline{z_{j}}A_{i}\right) \left(
z_{i}A_{j}^{\ast }-z_{j}A_{i}^{\ast }\right) \geq 0,
\end{equation*}%
from where results%
\begin{equation}
\left\vert z_{i}\right\vert ^{2}A_{j}A_{j}^{\ast }+\left\vert
z_{j}\right\vert ^{2}A_{i}A_{i}^{\ast }\geq \overline{z_{i}}%
z_{j}A_{j}A_{i}^{\ast }+\overline{z_{j}}z_{i}A_{i}A_{j}^{\ast }  \label{1.2}
\end{equation}%
for any $i,j\in \left\{ 1,\dots ,n\right\} .$

If we sum (\ref{1.2}) over $i$ from $1$ to $n$ we deduce%
\begin{equation}
\left( \sum_{i=1}^{n}\left\vert z_{i}\right\vert ^{2}\right)
A_{j}A_{j}^{\ast }+\left\vert z_{j}\right\vert ^{2}\left(
\sum_{i=1}^{n}A_{i}A_{i}^{\ast }\right) \geq z_{j}A_{j}\left( \sum_{i=1}^{n}%
\overline{z_{i}}A_{i}^{\ast }\right) +\left( \sum_{i=1}^{n}z_{i}A_{i}\right) 
\overline{z_{j}}A_{j}^{\ast },  \label{1.3}
\end{equation}%
for any $j\in \left\{ 1,\dots ,n\right\} .$

If we sum (\ref{1.3}) over $j$ from $1$ to $n,$ we deduce%
\begin{multline}
\sum_{i=1}^{n}\left\vert z_{i}\right\vert ^{2}\sum_{j=1}^{n}A_{j}A_{j}^{\ast
}+\sum_{j=1}^{n}\left\vert z_{j}\right\vert ^{2}\left(
\sum_{i=1}^{n}A_{i}A_{i}^{\ast }\right)  \label{1.4} \\
\geq \sum_{j=1}^{n}z_{j}A_{j}\sum_{i=1}^{n}\overline{z_{i}}A_{i}^{\ast
}+\left( \sum_{i=1}^{n}z_{i}A_{i}\right) \left( \sum_{j=1}^{n}\overline{z_{j}%
}A_{j}^{\ast }\right) ,
\end{multline}%
that is,%
\begin{equation}
\sum_{k=1}^{n}\left\vert z_{k}\right\vert ^{2}\sum_{k=1}^{n}A_{k}A_{k}^{\ast
}\geq \sum_{k=1}^{n}z_{k}A_{k}\sum_{k=1}^{n}\overline{z_{k}}A_{k}^{\ast
}=\left( \sum_{k=1}^{n}z_{k}A_{k}\right) \left(
\sum_{k=1}^{n}z_{k}A_{k}\right) ^{\ast }\geq 0,  \label{1.5}
\end{equation}%
and the theorem is proved.
\end{proof}

The following version of the Cauchy-Bunyakovsky-Schwarz inequality for norms
also holds \cite{SSD}.

\begin{corollary}
\label{c1.2}With the assumptions in Theorem \ref{t1.1}, one has%
\begin{equation}
\sum_{k=1}^{n}\left\vert z_{k}\right\vert ^{2}\left\Vert
\sum_{k=1}^{n}A_{k}A_{k}^{\ast }\right\Vert \geq \left\Vert
\sum_{k=1}^{n}z_{k}A_{k}\right\Vert ^{2}.  \label{1.6}
\end{equation}
\end{corollary}

\begin{proof}
The operators:%
\begin{equation*}
A:=\sum_{k=1}^{n}\left\vert z_{k}\right\vert
^{2}\sum_{k=1}^{n}A_{k}A_{k}^{\ast },\ \ B:=\left(
\sum_{k=1}^{n}z_{k}A_{k}\right) \left( \sum_{k=1}^{n}\overline{z_{k}}%
A_{k}\right) ^{\ast }
\end{equation*}%
are obviously self-adjoint, positive and by (\ref{1.1}), $A\geq B\geq 0.$
Thus $\left\Vert A\right\Vert \geq \left\Vert B\right\Vert $ and since,%
\begin{equation*}
\left\Vert A\right\Vert =\sum_{k=1}^{n}\left\vert z_{k}\right\vert
^{2}\left\Vert \sum_{k=1}^{n}A_{k}A_{k}^{\ast }\right\Vert 
\end{equation*}%
and%
\begin{equation*}
\left\Vert B\right\Vert =\left\Vert \sum_{k=1}^{n}z_{k}A_{k}\right\Vert ^{2}
\end{equation*}%
the corollary is proved.
\end{proof}

For other related results, see \cite{DM}.

The main aim of this paper is to point out other inequalities similar to (%
\ref{1.6}).

\section{Norm Inequalities}

The following result holds.

\begin{theorem}
\label{t2.1}Let $\alpha _{1},\dots ,\alpha _{n}\in \mathbb{K}$ and $%
A_{1},\dots ,A_{n}\in B\left( H\right) .$ Then one has the inequalities:%
\begin{multline}
\left\Vert \sum_{i=1}^{n}\alpha _{i}A_{i}\right\Vert ^{2}\leq \left\{ 
\begin{array}{l}
\max\limits_{i=\overline{1,n}}\left\vert \alpha _{i}\right\vert
^{2}\sum\limits_{i=1}^{n}\left\Vert A_{i}\right\Vert ^{2} \\ 
\\ 
\left( \sum\limits_{i=1}^{n}\left\vert \alpha _{i}\right\vert ^{2p}\right) ^{%
\frac{1}{p}}\left( \sum\limits_{i=1}^{n}\left\Vert A_{i}\right\Vert
^{2q}\right) ^{\frac{1}{q}}\text{ \ if \ }p>1,\ \frac{1}{p}+\frac{1}{q}=1;
\\ 
\\ 
\sum\limits_{i=1}^{n}\left\vert \alpha _{i}\right\vert ^{2}\max\limits_{i=%
\overline{1,n}}\left\Vert A_{i}\right\Vert ^{2}%
\end{array}%
\right.   \label{2.1} \\
+\left\{ 
\begin{array}{l}
\max\limits_{1\leq i\neq j\leq n}\left\{ \left\vert \alpha _{i}\right\vert
\left\vert \alpha _{j}\right\vert \right\} \sum\limits_{1\leq i\neq j\leq
n}\left\Vert A_{i}A_{j}^{\ast }\right\Vert  \\ 
\\ 
\left[ \left( \sum\limits_{i=1}^{n}\left\vert \alpha _{i}\right\vert
^{r}\right) ^{2}-\left( \sum\limits_{i=1}^{n}\left\vert \alpha
_{i}\right\vert ^{2r}\right) ^{\frac{1}{2}}\right] \left( \sum\limits_{1\leq
i\neq j\leq n}\left\Vert A_{i}A_{j}^{\ast }\right\Vert ^{s}\right) ^{\frac{1%
}{s}} \\ 
\text{\hfill if \ }r>1,\ \frac{1}{r}+\frac{1}{s}=1; \\ 
\left[ \left( \sum\limits_{i=1}^{n}\left\vert \alpha _{i}\right\vert \right)
^{2}-\left( \sum\limits_{i=1}^{n}\left\vert \alpha _{i}\right\vert
^{2}\right) \right] \max\limits_{1\leq i\neq j\leq n}\left\Vert
A_{i}A_{j}^{\ast }\right\Vert ,%
\end{array}%
\right. 
\end{multline}%
where (\ref{2.1}) should be seen as all the 9 possible configurations.
\end{theorem}

\begin{proof}
We have%
\begin{align}
0& \leq \left( \sum_{i=1}^{n}\alpha _{i}A_{i}\right) \left(
\sum_{i=1}^{n}\alpha _{i}A_{i}\right) ^{\ast }=\left( \sum_{i=1}^{n}\alpha
_{i}A_{i}\right) \left( \sum_{j=1}^{n}\overline{\alpha _{j}}A_{j}^{\ast
}\right)   \label{2.2} \\
& =\sum_{i=1}^{n}\sum_{j=1}^{n}\alpha _{i}\overline{\alpha _{j}}%
A_{i}A_{j}^{\ast }=\sum\limits_{i=1}^{n}\left\vert \alpha _{i}\right\vert
^{2}A_{i}A_{i}^{\ast }+\sum\limits_{1\leq i\neq j\leq n}\alpha _{i}\overline{%
\alpha _{j}}A_{i}A_{j}^{\ast }.  \notag
\end{align}%
Taking the norm in $\left( \ref{2.2}\right) $ and observing  that $%
\left\Vert UU^{\ast }\right\Vert =\left\Vert U\right\Vert ^{2}$ for any $%
U\in B\left( H\right) ,$ one has the inequality%
\begin{align}
\left\Vert \sum_{i=1}^{n}\alpha _{i}A_{i}\right\Vert ^{2}& =\left\Vert
\sum\limits_{i=1}^{n}\left\vert \alpha _{i}\right\vert ^{2}A_{i}A_{i}^{\ast
}+\sum\limits_{1\leq i\neq j\leq n}\alpha _{i}\overline{\alpha _{j}}%
A_{i}A_{j}^{\ast }\right\Vert   \label{2.3} \\
& \leq \sum\limits_{i=1}^{n}\left\vert \alpha _{i}\right\vert ^{2}\left\Vert
A_{i}A_{i}^{\ast }\right\Vert +\sum\limits_{1\leq i\neq j\leq n}\left\vert
\alpha _{i}\right\vert \left\vert \alpha _{j}\right\vert \left\Vert
A_{i}A_{j}^{\ast }\right\Vert   \notag \\
& =\sum\limits_{i=1}^{n}\left\vert \alpha _{i}\right\vert ^{2}\left\Vert
A_{i}\right\Vert ^{2}+\sum\limits_{1\leq i\neq j\leq n}\left\vert \alpha
_{i}\right\vert \left\vert \alpha _{j}\right\vert \left\Vert
A_{i}A_{j}^{\ast }\right\Vert .  \notag
\end{align}%
Using H\"{o}lder's inequality, we may write that:%
\begin{equation}
\sum\limits_{i=1}^{n}\left\vert \alpha _{i}\right\vert ^{2}\left\Vert
A_{i}\right\Vert ^{2}\leq \left\{ 
\begin{array}{l}
\max\limits_{i=\overline{1,n}}\left\vert \alpha _{i}\right\vert
^{2}\sum\limits_{1\leq i\neq j\leq n}\left\Vert A_{i}\right\Vert ^{2} \\ 
\\ 
\left( \sum\limits_{i=1}^{n}\left\vert \alpha _{i}\right\vert ^{2p}\right) ^{%
\frac{1}{p}}\left( \sum\limits_{i=1}^{n}\left\Vert A_{i}\right\Vert
^{2q}\right) ^{\frac{1}{q}}\text{ \ if \ }p>1,\ \frac{1}{p}+\frac{1}{q}=1;
\\ 
\\ 
\sum\limits_{i=1}^{n}\left\vert \alpha _{i}\right\vert ^{2}\max\limits_{i=%
\overline{1,n}}\left\Vert A_{i}\right\Vert ^{2}.%
\end{array}%
\right.   \label{2.4}
\end{equation}%
Also, H\"{o}lder's inequality for double sums produces%
\begin{align}
\sum\limits_{1\leq i\neq j\leq n}\left\vert \alpha _{i}\right\vert
\left\vert \alpha _{j}\right\vert \left\Vert A_{i}A_{j}^{\ast }\right\Vert &
\leq \left\{ 
\begin{array}{l}
\max\limits_{1\leq i\neq j\leq n}\left\{ \left\vert \alpha _{i}\right\vert
\left\vert \alpha _{j}\right\vert \right\} \sum\limits_{1\leq i\neq j\leq
n}\left\Vert A_{i}A_{j}^{\ast }\right\Vert  \\ 
\\ 
\left( \sum\limits_{1\leq i\neq j\leq n}\left\vert \alpha _{i}\right\vert
^{r}\left\vert \alpha _{j}\right\vert ^{r}\right) ^{\frac{1}{r}}\left(
\sum\limits_{1\leq i\neq j\leq n}\left\Vert A_{i}A_{j}^{\ast }\right\Vert
^{s}\right) ^{\frac{1}{s}} \\ 
\text{\hfill if \ }r>1,\ \frac{1}{r}+\frac{1}{s}=1; \\ 
\sum\limits_{1\leq i\neq j\leq n}\left\vert \alpha _{i}\right\vert
\left\vert \alpha _{j}\right\vert \max\limits_{1\leq i\neq j\leq
n}\left\Vert A_{i}A_{j}^{\ast }\right\Vert ,%
\end{array}%
\right.   \label{2.5} \\
& =\left\{ 
\begin{array}{l}
\max\limits_{1\leq i\neq j\leq n}\left\{ \left\vert \alpha _{i}\right\vert
\left\vert \alpha _{j}\right\vert \right\} \sum\limits_{1\leq i\neq j\leq
n}\left\Vert A_{i}A_{j}^{\ast }\right\Vert  \\ 
\\ 
\left[ \left( \sum\limits_{i=1}^{n}\left\vert \alpha _{i}\right\vert
^{r}\right) ^{2}-\left( \sum\limits_{i=1}^{n}\left\vert \alpha
_{i}\right\vert ^{2r}\right) ^{\frac{1}{2}}\right] \left(
\sum\limits_{i=1}^{n}\left\Vert A_{i}A_{j}^{\ast }\right\Vert ^{s}\right) ^{%
\frac{1}{s}} \\ 
\text{\hfill if \ }r>1,\ \frac{1}{r}+\frac{1}{s}=1; \\ 
\left[ \left( \sum\limits_{i=1}^{n}\left\vert \alpha _{i}\right\vert \right)
^{2}-\left( \sum\limits_{i=1}^{n}\left\vert \alpha _{i}\right\vert
^{2}\right) \right] \max\limits_{1\leq i\neq j\leq n}\left\Vert
A_{i}A_{j}^{\ast }\right\Vert ,%
\end{array}%
\right.   \notag
\end{align}%
Using (\ref{2.3}) and (\ref{2.4}), (\ref{2.5}) one deduces the desired
inequality (\ref{2.1}).
\end{proof}

The following corollaries are natural consequences.

\begin{corollary}
\label{c2.2}With the assumptions of Theorem \ref{t2.1}, one has the
inequality%
\begin{equation}
\left\Vert \sum_{i=1}^{n}\alpha _{i}A_{i}\right\Vert \leq \max\limits_{i=%
\overline{1,n}}\left\vert \alpha _{i}\right\vert \left(
\sum\limits_{i,j=1}^{n}\left\Vert A_{i}A_{j}^{\ast }\right\Vert \right) ^{%
\frac{1}{2}}.  \label{2.6}
\end{equation}
\end{corollary}

\begin{proof}
Follows by the first line in (\ref{2.1}) on taking into account that 
\begin{equation*}
\max\limits_{1\leq i\neq j\leq n}\left\{ \left\vert \alpha _{i}\right\vert
\left\vert \alpha _{j}\right\vert \right\} \leq \max_{i=\overline{1,n}%
}\left\vert \alpha _{i}\right\vert ^{2},
\end{equation*}%
and 
\begin{equation*}
\sum\limits_{i,j=1}^{n}\left\Vert A_{i}A_{j}^{\ast }\right\Vert
=\sum_{i=1}^{n}\left\Vert A_{i}\right\Vert ^{2}+\sum\limits_{1\leq i\neq
j\leq n}\left\Vert A_{i}A_{j}^{\ast }\right\Vert .
\end{equation*}
\end{proof}

\begin{corollary}
\label{c2.3}With the assumptions in Theorem \ref{t2.1}, one has the
inequality:%
\begin{multline}
\left\Vert \sum_{i=1}^{n}\alpha _{i}A_{i}\right\Vert ^{2}  \label{2.7} \\
\leq \left( \sum\limits_{i=1}^{n}\left\vert \alpha _{i}\right\vert
^{2p}\right) ^{\frac{1}{p}}\left[ \left( \sum_{i=1}^{n}\left\Vert
A_{i}\right\Vert ^{2q}\right) ^{\frac{1}{q}}+\left( n-1\right) \left(
\sum\limits_{1\leq i\neq j\leq n}\left\Vert A_{i}A_{j}^{\ast }\right\Vert
^{q}\right) ^{\frac{1}{q}}\right] ,
\end{multline}%
where $p>1,$ $\frac{1}{p}+\frac{1}{q}=1.$
\end{corollary}

\begin{proof}
Using the Cauchy-Bunyakovsky-Schwarz inequality for positive numbers%
\begin{equation*}
\left( \sum\limits_{i=1}^{n}a_{i}\right) ^{2}\leq
n\sum\limits_{i=1}^{n}a_{i}^{2}
\end{equation*}%
we may write that%
\begin{align*}
\left( \sum\limits_{i=1}^{n}\left\vert \alpha _{i}\right\vert ^{p}\right)
^{2}-\sum\limits_{i=1}^{n}\left\vert \alpha _{i}\right\vert ^{2p}& \leq
n\sum\limits_{i=1}^{n}\left\vert \alpha _{i}\right\vert
^{2p}-\sum\limits_{i=1}^{n}\left\vert \alpha _{i}\right\vert ^{2p} \\
& =\left( n-1\right) \sum\limits_{i=1}^{n}\left\vert \alpha _{i}\right\vert
^{2p}.
\end{align*}%
Now, using the second line in (\ref{2.1}) for $r=p,$ $s=q,$ we deduce the
desired result (\ref{2.7}).
\end{proof}

\begin{corollary}
\label{c2.4}With the assumptions in Theorem \ref{t2.1}, one has the
inequality%
\begin{equation}
\left\Vert \sum_{i=1}^{n}\alpha _{i}A_{i}\right\Vert ^{2}\leq
\sum\limits_{i=1}^{n}\left\vert \alpha _{i}\right\vert ^{2}\left[ \max_{i=%
\overline{1,n}}\left\Vert A_{i}\right\Vert ^{2}+\left( n-1\right)
\max\limits_{1\leq i\neq j\leq n}\left\Vert A_{i}A_{j}^{\ast }\right\Vert %
\right] .  \label{2.8}
\end{equation}
\end{corollary}

\begin{proof}
Follows by the third line of (\ref{2.1}) on taking into account that 
\begin{equation*}
\left[ \left( \sum\limits_{i=1}^{n}\left\vert \alpha _{i}\right\vert \right)
^{2}-\sum\limits_{i=1}^{n}\left\vert \alpha _{i}\right\vert ^{2}\right] ^{%
\frac{1}{2}}\leq \left( n-1\right) \sum\limits_{i=1}^{n}\left\vert \alpha
_{i}\right\vert ^{2}.
\end{equation*}
\end{proof}

Another interesting particular case is embodied in the following corollary
as well.

\begin{corollary}
\label{c2.5}With the assumptions in Theorem \ref{t2.1}, one has the
inequality 
\begin{equation}
\left\Vert \sum_{i=1}^{n}\alpha _{i}A_{i}\right\Vert ^{2}\leq
\sum\limits_{i=1}^{n}\left\vert \alpha _{i}\right\vert ^{2}\left[ \max_{i=%
\overline{1,n}}\left\Vert A_{i}\right\Vert ^{2}+\left( \sum\limits_{1\leq
i\neq j\leq n}\left\Vert A_{i}A_{j}^{\ast }\right\Vert ^{2}\right) ^{\frac{1%
}{2}}\right] .  \label{2.9}
\end{equation}
\end{corollary}

\begin{proof}
It is obvious that%
\begin{equation*}
\left[ \left( \sum\limits_{i=1}^{n}\left\vert \alpha _{i}\right\vert
^{2}\right) ^{2}-\sum\limits_{i=1}^{n}\left\vert \alpha _{i}\right\vert ^{4}%
\right] ^{\frac{1}{2}}\leq \sum\limits_{i=1}^{n}\left\vert \alpha
_{i}\right\vert ^{2}.
\end{equation*}%
Thus, combining the third line in the first bracket in (\ref{2.1}) with the
second line for $r=s=2$ in the second bracket, the inequality (\ref{2.9}) is
obtained.
\end{proof}

\begin{remark}
\label{r2.6}If one is interested in obtaining bounds in terms of $%
\sum_{i=1}^{n}\left\vert \alpha _{i}\right\vert ^{2},$ there are other
possibilities as shown below. Obviously, since%
\begin{equation*}
\max\limits_{1\leq i\neq j\leq n}\left\{ \left\vert \alpha _{i}\right\vert
\left\vert \alpha _{j}\right\vert \right\} \leq \max_{i=\overline{1,n}%
}\left\vert \alpha _{i}\right\vert ^{2}\leq \sum\limits_{i=1}^{n}\left\vert
\alpha _{i}\right\vert ^{2}.
\end{equation*}%
then, by (\ref{2.1}), in choosing the third line in the first bracket with
the first line in the second bracket, one would obtain%
\begin{equation}
\left\Vert \sum_{i=1}^{n}\alpha _{i}A_{i}\right\Vert ^{2}\leq
\sum\limits_{i=1}^{n}\left\vert \alpha _{i}\right\vert ^{2}\left[ \max_{i=%
\overline{1,n}}\left\Vert A_{i}\right\Vert ^{2}+\sum\limits_{1\leq i\neq
j\leq n}\left\Vert A_{i}A_{j}^{\ast }\right\Vert \right] .  \label{2.10}
\end{equation}%
Also, it is evident that%
\begin{equation*}
\left[ \left( \sum\limits_{i=1}^{n}\left\vert \alpha _{i}\right\vert
^{r}\right) ^{2}-\left( \sum\limits_{i=1}^{n}\left\vert \alpha
_{i}\right\vert ^{2r}\right) ^{\frac{1}{2}}\right] ^{\frac{1}{r}}\leq \left(
\sum\limits_{i=1}^{n}\left\vert \alpha _{i}\right\vert ^{r}\right) ^{\frac{2%
}{r}}.
\end{equation*}%
By the monotonicity of the power mean $\left( \frac{1}{n}%
\sum_{i=1}^{n}a_{i}^{m}\right) ^{\frac{1}{m}}$ as a function of $m\in 
\mathbb{R}$, we have%
\begin{equation*}
\left( \frac{\sum_{i=1}^{n}\left\vert \alpha _{i}\right\vert ^{r}}{n}\right)
^{\frac{1}{r}}\leq \left( \frac{\sum_{i=1}^{n}\left\vert \alpha
_{i}\right\vert ^{2}}{n}\right) ^{\frac{1}{2}},\ \ \ \ 1<r\leq 2,
\end{equation*}%
giving%
\begin{equation*}
\left( \sum\limits_{i=1}^{n}\left\vert \alpha _{i}\right\vert ^{r}\right) ^{%
\frac{2}{r}}\leq n^{\frac{2}{r}-1}\sum\limits_{i=1}^{n}\left\vert \alpha
_{i}\right\vert ^{2}.
\end{equation*}%
Thus, using the third line in the first bracket of (\ref{2.1}) combined with
the second line in the second bracket for $1\leq r\leq 2$ , $\frac{1}{s}+%
\frac{1}{r}=1,$ we deduce%
\begin{equation}
\left\Vert \sum_{i=1}^{n}\alpha _{i}A_{i}\right\Vert ^{2}\leq
\sum\limits_{i=1}^{n}\left\vert \alpha _{i}\right\vert ^{2}\left[ \max_{i=%
\overline{1,n}}\left\Vert A_{i}\right\Vert ^{2}+n^{\frac{2}{r}-1}\left(
\sum\limits_{1\leq i\neq j\leq n}\left\Vert A_{i}A_{j}^{\ast }\right\Vert
^{s}\right) ^{\frac{1}{s}}\right] .  \label{2.11}
\end{equation}%
Note that for $r=s=2,$ we recapture (\ref{2.9}).
\end{remark}

The following particular result also holds.

\begin{proposition}
\label{p2.6}Let $\alpha _{1},\dots ,\alpha _{n}\in \mathbb{K}$ and $%
A_{1},\dots ,A_{n}\in B\left( H\right) $ with the property that $%
A_{i}A_{j}^{\ast }=0$ for any $i\neq j,$ $i,j\in \left\{ 1,\dots ,n\right\}
. $ Then one has the inequality;%
\begin{equation}
\left\Vert \sum_{i=1}^{n}\alpha _{i}A_{i}\right\Vert \leq \left\{ 
\begin{array}{l}
\max\limits_{i=\overline{1,n}}\left\vert \alpha _{i}\right\vert \left(
\sum\limits_{i=1}^{n}\left\Vert A_{i}\right\Vert ^{2}\right) ^{\frac{1}{2}},
\\ 
\\ 
\left( \sum\limits_{i=1}^{n}\left\vert \alpha _{i}\right\vert ^{2p}\right) ^{%
\frac{1}{2p}}\left( \sum\limits_{i=1}^{n}\left\Vert A_{i}\right\Vert
^{2q}\right) ^{\frac{1}{2q}}\text{ \ if \ }p>1,\ \frac{1}{p}+\frac{1}{q}=1;
\\ 
\\ 
\left( \sum\limits_{i=1}^{n}\left\vert \alpha _{i}\right\vert ^{2}\right) ^{%
\frac{1}{2}}\max\limits_{i=\overline{1,n}}\left\Vert A_{i}\right\Vert .%
\end{array}%
\right.  \label{2.12}
\end{equation}
\end{proposition}

If by $M\left( \mathbf{\alpha },\mathbf{A}\right) $ we denote any of the
bounds provided by (\ref{2.1}), (\ref{2.6}), (\ref{2.7}), (\ref{2.8}), (\ref%
{2.9}), (\ref{2.10}) or (\ref{2.11}), then we may state the following
proposition as well.

\begin{proposition}
\label{p2.7}Under the assumptions of Theorem \ref{t2.1}, we have:

\begin{enumerate}
\item[(i)] For any $x\in H$%
\begin{equation}
\left\Vert \sum_{i=1}^{n}\alpha _{i}A_{i}x\right\Vert ^{2}\leq \left\Vert
x\right\Vert ^{2}M\left( \mathbf{\alpha },\mathbf{A}\right) .  \label{2.13}
\end{equation}

\item[(ii)] For any $x,y\in H,$%
\begin{equation}
\left\vert \sum_{i=1}^{n}\alpha _{i}\left\langle A_{i}x,y\right\rangle
\right\vert ^{2}\leq \left\Vert x\right\Vert ^{2}\left\Vert y\right\Vert
^{2}M\left( \mathbf{\alpha },\mathbf{A}\right) .  \label{2.14}
\end{equation}
\end{enumerate}
\end{proposition}

\begin{proof}

\begin{enumerate}
\item[(i)] Obviously, 
\begin{align*}
\left\Vert \sum_{i=1}^{n}\alpha _{i}A_{i}x\right\Vert ^{2}& =\left\Vert
\left( \sum_{i=1}^{n}\alpha _{i}A_{i}\right) \left( x\right) \right\Vert
^{2}\leq \left\Vert \sum_{i=1}^{n}\alpha _{i}A_{i}\right\Vert ^{2}\left\Vert
x\right\Vert ^{2} \\
& \leq M\left( \mathbf{\alpha },\mathbf{A}\right) \left\Vert x\right\Vert
^{2}.
\end{align*}

\item[(ii)] We have%
\begin{equation*}
\left\vert \sum_{i=1}^{n}\alpha _{i}\left\langle A_{i}x,y\right\rangle
\right\vert ^{2}=\left\vert \left\langle \sum_{i=1}^{n}\alpha
_{i}A_{i}x,y\right\rangle \right\vert ^{2}=\left\Vert \sum_{i=1}^{n}\alpha
_{i}A_{i}x\right\Vert ^{2}\left\Vert y\right\Vert ^{2},
\end{equation*}%
which, by (i), gives the desired result (\ref{2.14}).
\end{enumerate}
\end{proof}

\section{Inequalities for Vectors in Hilbert Spaces}

We consider the non zero vectors $y_{1},\dots ,y_{n}\in H.$ Define the
operators%
\begin{equation*}
A_{i}:H\rightarrow H,\ \ A_{i}x=\frac{\left( x,y_{i}\right) }{\left\Vert
y_{i}\right\Vert }\cdot y_{i},\ \ \ i\in \left\{ 1,\dots ,n\right\} .
\end{equation*}%
Since%
\begin{equation}
\left\Vert A_{i}\right\Vert =\sup_{\left\Vert x\right\Vert =1}\left\Vert
A_{i}x\right\Vert =\sup_{\left\Vert x\right\Vert =1}\left\vert \left(
x,y_{i}\right) \right\vert =\left\Vert y_{i}\right\Vert ,\ \ \ i\in \left\{
1,\dots ,n\right\}   \label{3.1}
\end{equation}%
then $A_{i}$ are bounded linear operators in $H.$ Also, since%
\begin{equation}
\left( A_{i}x,x\right) =\left( \frac{\left( x,y_{i}\right) y_{i}}{\left\Vert
y_{i}\right\Vert },x\right) =\frac{\left\vert \left( x,y_{i}\right)
\right\vert ^{2}}{\left\Vert y_{i}\right\Vert }\geq 0,\ \ x\in H,\ \ i\in
\left\{ 1,\dots ,n\right\}   \label{3.2}
\end{equation}%
and%
\begin{align*}
\left( A_{i}x,z\right) & =\left( \frac{\left( x,y_{i}\right) y_{i}}{%
\left\Vert y_{i}\right\Vert },z\right) =\frac{\left( x,y_{i}\right) \left(
y_{i},z\right) }{\left\Vert y_{i}\right\Vert }, \\
\left( x,A_{i}z\right) & =\left( x,\frac{\left( z,y_{i}\right) y_{i}}{%
\left\Vert y_{i}\right\Vert }\right) =\frac{\left( x,y_{i}\right) \overline{%
\left( z,y_{i}\right) }}{\left\Vert y_{i}\right\Vert }=\frac{\left(
x,y_{i}\right) \left( y_{i},z\right) }{\left\Vert y_{i}\right\Vert },
\end{align*}%
giving%
\begin{equation}
\left( A_{i}x,z\right) =\left( x,A_{i}z\right) ,\ \ \ x,z\in H,\ \ \ i\in
\left\{ 1,\dots ,n\right\} ,  \label{3.3}
\end{equation}%
we may conclude that $A_{i}$ $\left( i=1,\dots ,n\right) $ are positive
self-adjoint operators on $H.$

Since, for any $x\in H,$ one has%
\begin{align*}
\left\Vert \left( A_{i}A_{j}\right) \left( x\right) \right\Vert &
=\left\Vert \left( A_{i}\right) \left( A_{j}x\right) \right\Vert =\left\Vert
A_{i}\left( \frac{\left( x,y_{j}\right) y_{j}}{\left\Vert y_{j}\right\Vert }%
\right) \right\Vert \\
& =\frac{\left\vert \left( x,y_{j}\right) \right\vert }{\left\Vert
y_{j}\right\Vert }\left\Vert A_{i}y_{j}\right\Vert =\frac{\left\vert \left(
x,y_{j}\right) \right\vert }{\left\Vert y_{j}\right\Vert }\cdot \frac{%
\left\vert \left( y_{j},y_{i}\right) \right\vert \left\Vert y_{j}\right\Vert 
}{\left\Vert y_{i}\right\Vert } \\
& =\frac{\left\vert \left( x,y_{j}\right) \right\vert \left\vert \left(
y_{j},y_{i}\right) \right\vert }{\left\Vert y_{i}\right\Vert },\ \ i,j\in
\left\{ 1,\dots ,n\right\} ,
\end{align*}%
we deduce that%
\begin{equation}
\left\Vert A_{i}A_{j}\right\Vert =\sup_{\left\Vert \alpha \right\Vert =1}%
\frac{\left\vert \left( x,y_{j}\right) \right\vert \left\vert \left(
y_{j},y_{i}\right) \right\vert }{\left\Vert y_{i}\right\Vert }=\left\vert
\left( y_{i},y_{j}\right) \right\vert ;\ \ i,j\in \left\{ 1,\dots ,n\right\}
.  \label{3.4}
\end{equation}%
If $\left( y_{i}\right) _{i=\overline{1,n}}$ is an orthogonal family on $H,$
then $\left\Vert A_{i}\right\Vert =1$ and $A_{i}A_{j}=0$ for $i,j\in \left\{
1,\dots ,n\right\} ,$ $i\neq j.$

The following inequality for vectors holds.

\begin{theorem}
\label{t3.1}Let $x,y_{1},\dots ,y_{n}\in H$ and $\alpha _{n},\dots ,\alpha
_{n}\in \mathbb{K}$. Then one has the inequalities:%
\begin{multline}
\left\Vert \sum_{i=1}^{n}\alpha _{i}\frac{\left( x,y_{i}\right) }{\left\Vert
y_{i}\right\Vert }y_{i}\right\Vert ^{2}  \label{3.5} \\
\leq \left\Vert x\right\Vert ^{2}\times \left\{ 
\begin{array}{l}
\max\limits_{i=\overline{1,n}}\left\vert \alpha _{i}\right\vert
^{2}\sum\limits_{i=1}^{n}\left\Vert y_{i}\right\Vert ^{2} \\ 
\\ 
\left( \sum\limits_{i=1}^{n}\left\vert \alpha _{i}\right\vert ^{2p}\right) ^{%
\frac{1}{p}}\left( \sum\limits_{i=1}^{n}\left\Vert y_{i}\right\Vert
^{2q}\right) ^{\frac{1}{q}}\text{ \ if \ }p>1,\ \frac{1}{p}+\frac{1}{q}=1;
\\ 
\\ 
\sum\limits_{i=1}^{n}\left\vert \alpha _{i}\right\vert ^{2}\max\limits_{i=%
\overline{1,n}}\left\Vert y_{i}\right\Vert ^{2}%
\end{array}%
\right.  \\
+\left\Vert x\right\Vert ^{2}\times \left\{ 
\begin{array}{l}
\max\limits_{1\leq i\neq j\leq n}\left\{ \left\vert \alpha _{i}\right\vert
\left\vert \alpha _{j}\right\vert \right\} \sum\limits_{1\leq i\neq j\leq
n}\left\vert \left( y_{i},y_{j}\right) \right\vert  \\ 
\\ 
\left[ \left( \sum\limits_{i=1}^{n}\left\vert \alpha _{i}\right\vert
^{r}\right) ^{2}-\sum\limits_{i=1}^{n}\left\vert \alpha _{i}\right\vert ^{2r}%
\right] ^{r}\left( \sum\limits_{1\leq i\neq j\leq n}\left\vert \left(
y_{i},y_{j}\right) \right\vert ^{s}\right) ^{\frac{1}{s}} \\ 
\text{\hfill if \ }r>1,\ \frac{1}{r}+\frac{1}{s}=1; \\ 
\left[ \left( \sum\limits_{i=1}^{n}\left\vert \alpha _{i}\right\vert \right)
^{2}-\sum\limits_{i=1}^{n}\left\vert \alpha _{i}\right\vert ^{2}\right]
\max\limits_{1\leq i\neq j\leq n}\left\vert \left( y_{i},y_{j}\right)
\right\vert .%
\end{array}%
\right. 
\end{multline}
\end{theorem}

\begin{proof}
Follows by Theorem \ref{t2.1} and Proposition \ref{p2.7}, (i) on choosing $%
A_{i}=\frac{\left( \cdot ,y_{i}\right) }{\left\Vert y_{i}\right\Vert }y_{i}$
and taking into account that $\left\Vert A_{i}\right\Vert =\left\Vert
y_{i}\right\Vert ,$ 
\begin{equation*}
\left\Vert A_{i}A_{j}^{\ast }\right\Vert =\left\vert \left(
y_{i},y_{j}\right) \right\vert ,\ \ \ \ \ i,j\in \left\{ 1,\dots ,n\right\} .
\end{equation*}%
We omit the details.
\end{proof}

Using Corollaries \ref{c2.2}--\ref{c2.5} and Remark \ref{r2.6}, we may state
the following particular inequalities:%
\begin{equation}
\left\Vert \sum_{i=1}^{n}\alpha _{i}\frac{\left( x,y_{i}\right) }{\left\Vert
y_{i}\right\Vert }y_{i}\right\Vert \leq \left\Vert x\right\Vert
\max\limits_{i=\overline{1,n}}\left\vert \alpha _{i}\right\vert \left(
\sum_{i,j=1}^{n}\left\vert \left( y_{i},y_{j}\right) \right\vert \right) ^{%
\frac{1}{2}};  \label{3.6}
\end{equation}%
\begin{multline}
\left\Vert \sum_{i=1}^{n}\alpha _{i}\frac{\left( x,y_{i}\right) }{\left\Vert
y_{i}\right\Vert }y_{i}\right\Vert ^{2}  \label{3.7} \\
\leq \left\Vert x\right\Vert ^{2}\left[ \left(
\sum\limits_{i=1}^{n}\left\vert \alpha _{i}\right\vert ^{2p}\right) ^{\frac{1%
}{p}}\left[ \left( \sum\limits_{i=1}^{n}\left\Vert y_{i}\right\Vert
^{2q}\right) ^{\frac{1}{q}}+\left( n-1\right) \left( \sum\limits_{1\leq
i\neq j\leq n}\left\vert \left( y_{i},y_{j}\right) \right\vert ^{q}\right) ^{%
\frac{1}{q}}\right] \right] ,
\end{multline}%
where $p>1,$ $\frac{1}{p}+\frac{1}{q}=1;$%
\begin{equation}
\left\Vert \sum_{i=1}^{n}\alpha _{i}\frac{\left( x,y_{i}\right) }{\left\Vert
y_{i}\right\Vert }y_{i}\right\Vert ^{2}\leq \left\Vert x\right\Vert
^{2}\sum\limits_{i=1}^{n}\left\vert \alpha _{i}\right\vert ^{2}\left[
\max\limits_{i=\overline{1,n}}\left\Vert y_{i}\right\Vert ^{2}+\left(
n-1\right) \max\limits_{1\leq i\neq j\leq n}\left\vert \left(
y_{i},y_{j}\right) \right\vert \right] ;  \label{3.8}
\end{equation}%
\begin{equation}
\left\Vert \sum_{i=1}^{n}\alpha _{i}\frac{\left( x,y_{i}\right) }{\left\Vert
y_{i}\right\Vert }y_{i}\right\Vert ^{2}\leq \left\Vert x\right\Vert
^{2}\sum\limits_{i=1}^{n}\left\vert \alpha _{i}\right\vert ^{2}\left[
\max\limits_{i=\overline{1,n}}\left\Vert y_{i}\right\Vert ^{2}+\left(
\sum\limits_{1\leq i\neq j\leq n}\left\vert \left( y_{i},y_{j}\right)
\right\vert ^{2}\right) ^{\frac{1}{2}}\right] ;  \label{3.9}
\end{equation}%
\begin{equation}
\left\Vert \sum_{i=1}^{n}\alpha _{i}\frac{\left( x,y_{i}\right) }{\left\Vert
y_{i}\right\Vert }y_{i}\right\Vert ^{2}\leq \left\Vert x\right\Vert
^{2}\sum\limits_{i=1}^{n}\left\vert \alpha _{i}\right\vert ^{2}\left[
\max\limits_{i=\overline{1,n}}\left\Vert y_{i}\right\Vert
^{2}+\sum\limits_{1\leq i\neq j\leq n}\left\vert \left( y_{i},y_{j}\right)
\right\vert \right] ;  \label{3.10}
\end{equation}%
\begin{multline}
\left\Vert \sum_{i=1}^{n}\alpha _{i}\frac{\left( x,y_{i}\right) }{\left\Vert
y_{i}\right\Vert }y_{i}\right\Vert ^{2}  \label{3.11} \\
\leq \left\Vert x\right\Vert ^{2}\sum\limits_{i=1}^{n}\left\vert \alpha
_{i}\right\vert ^{2}\left[ \max\limits_{i=\overline{1,n}}\left\Vert
y_{i}\right\Vert ^{2}+n^{\frac{2}{r}-1}\left( \sum\limits_{1\leq i\neq j\leq
n}\left\vert \left( y_{i},y_{j}\right) \right\vert ^{s}\right) ^{\frac{1}{s}}%
\right] ,
\end{multline}%
where $1<r\leq 2,$ $\frac{1}{s}+\frac{1}{r}=1.$

\begin{remark}
The choice $\alpha _{i}=\left\Vert y_{i}\right\Vert $ $\left( i=1,\dots
,n\right) $ will produce some interesting bounds for%
\begin{equation*}
\left\Vert \sum_{i=1}^{n}\left( x,y_{i}\right) y_{i}\right\Vert ^{2}.
\end{equation*}%
We omit the details.
\end{remark}

\end{document}